\documentclass[nomath,noamsfonts]{conm-p-l}
\usepackage{graphicx}

\newcommand\ds{\displaystyle}
\newcommand\bs{$\backslash$}
\newcommand\Disc{\scalebox{0.5}{\includegraphics{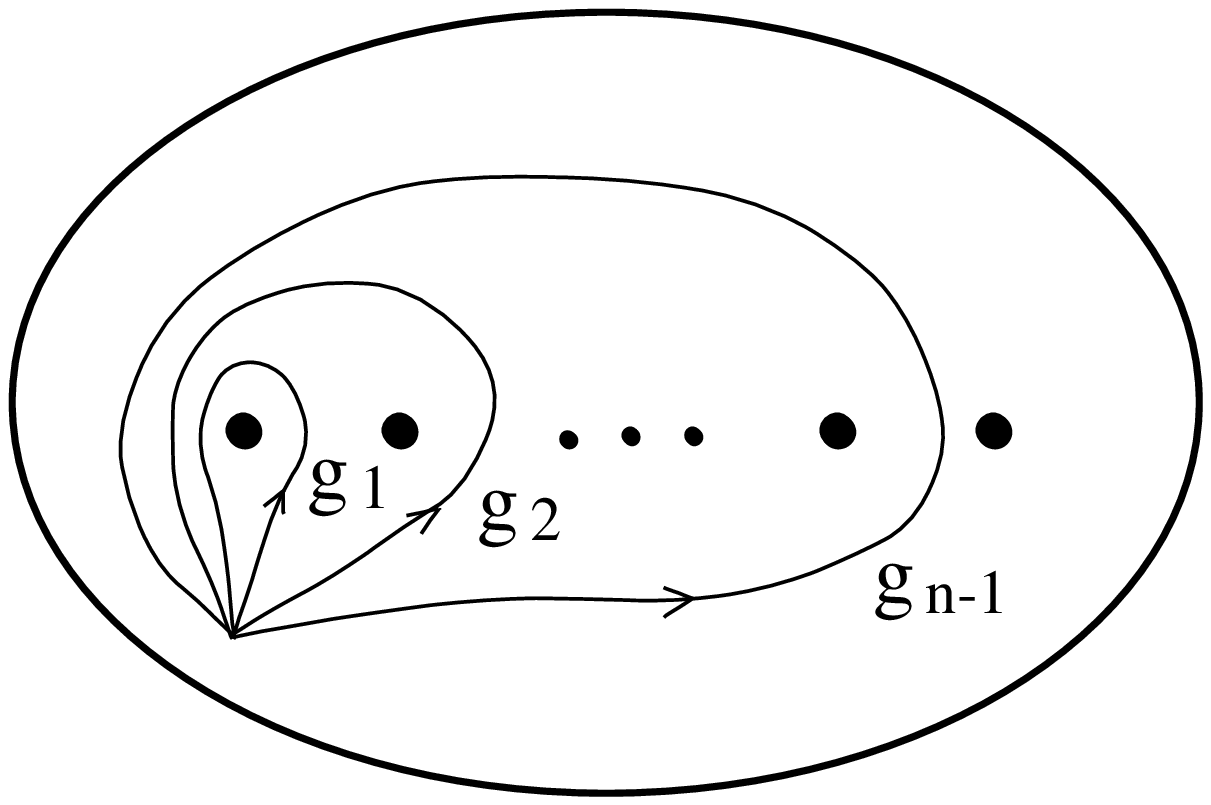}}}
\newcommand\Braid{\scalebox{0.5}{\includegraphics{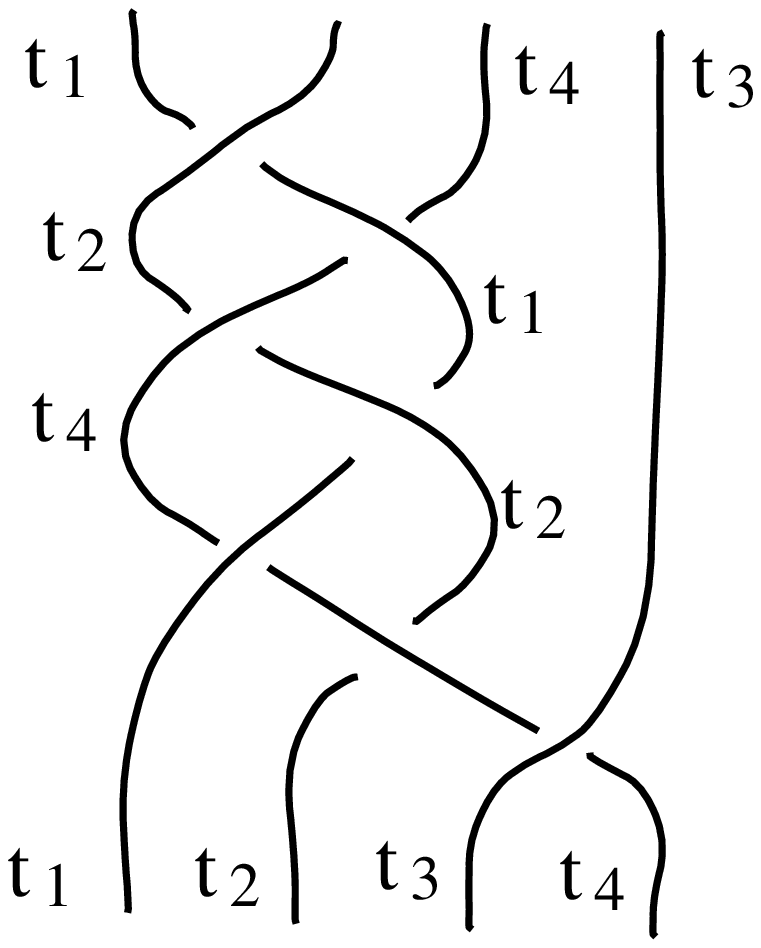}}}

\newtheorem{theorem}{Theorem}

\begin{document}


\title{The Multivariable Alexander
Polynomial for a Closed Braid}
\author{H. R. Morton}
\address{ Department of Mathematical Sciences\\
 University of Liverpool\\
 Liverpool L69 3BX\\ England.}
\email{morton@liv.ac.uk}
\urladdr{http://www.liv.ac.uk/\~{ }su14}
 \thanks{This paper was prepared during the low-dimensional
topology program at MSRI, Berkeley in 1996-97. I am grateful to MSRI for their
support.}
\begin{abstract}
A simple multivariable version of the reduced Burau matrix  is
constructed for any braid. It is shown how the multivariable Alexander
polynomial for the closure of the braid can be found directly from this matrix.
\end{abstract} 

\maketitle

\section{Introduction}

It has been known for some time that the Alexander polynomial of a closed braid
$\hat\beta$ can be found from the Burau matrix $\beta$, \cite{Birman}. This
relation was extended in \cite{Morton} to present the $2$-variable
Alexander polynomial $\Delta_{\hat{\beta}\cup A}(t,x)$ of the link consisting
of the closed braid $\hat\beta$ together with its axis $A$  as the
characteristic polynomial, $\det(I-x\overline{B}_\beta(t))$, of the reduced
$(n-1)\times(n-1)$ Burau matrix $\overline{B}_\beta(t)$ of the braid $\beta$.
The Alexander polynomial of $\hat\beta$ can be recovered by  applying
the Torres-Fox formula to the $2$-variable polynomial to get the equation
\[{\Delta_{\hat{\beta}}(t)\over 1-t}= {\Delta_{\hat{\beta}\cup A}(t,1)\over
1-t^n}.\]

In this paper I give a similar method for finding the multivariable Alexander
polynomial of a link $L$ presented as the closure of a braid $\beta$. The main
ingredient is a readily constructed multivariable version of the reduced Burau
matrices.  Other versions of `coloured' Burau matrices have been developed, for
example by Penne, \cite{Penne}, which can be interpreted as determining linear
presentations of suitably extended versions of the braid group.

  The most useful feature of the matrices which are used here is that they are
extremely simple to remember and they give an immediate and very
straightforward construction of the Alexander polynomial of a closed braid and
axis, leading at once to the polynomial of the closed braid. For a pure braid
the resulting  matrix is conjugate to a reduced version of the Gassner matrix;
the construction given here has the advantage that it applies to any  braid
which presents the link, and does not require the braid to be rewritten in any
special form.

An implementation of this calculation by a Maple procedure, which  returns
the multivariable Alexander polynomial of $\hat\beta$ given the braid $\beta$,
was made in early 1996 by a Liverpool MSc student, Julian Hodgson, and is
usually available from the Liverpool knot-theory website,{\tt \bs
http:\bs\bs  www.liv.ac.uk\bs PureMaths\bs knots.html}.

\section{The multivariable Alexander polynomial} Given a
homomorphism $\varphi:G_L\to H$ from the group $G_L$ of an oriented link $L$ to
an abelian group $H$ we can use Fox's free differential calculus on a
presentation of  $G_L$ to find an invariant of $L$ which lies in the
group ring
${\bf Z}[H]$.
 The full multivariable Alexander polynomial $\Delta_L$ arises  when $H$ is the
abelianisation $G_L/G'_L$ and $\varphi$ is the natural projection. In this case
$H$ is the free abelian group on $k$ generators, $t_1,\ldots,t_k$ say, where $L$
has $k$ components $L_1,\ldots,L_k$, and $\varphi(x_i)=t_i$ for every oriented
meridian $x_i$ of the component $L_i$. Then  $\Delta_L\in {\bf Z}[G_L/G'_L]$
is  a Laurent polynomial in $t_1,\ldots,t_k$. Any other homomorphism 
$\varphi:G_L\to H$  factors through a homomorphism
$\overline{\varphi}:G_L/G'_L\to H$
 and the resulting invariant is given by substituting the 
images of
$\overline{\varphi}(t_i)$ in $\Delta_L$, \cite{CF}.

\subsection{ A coloured Burau matrix for $\beta$.}

Label the individual strings of $\beta\in B_n$ by $t_1,\ldots,t_n$, putting the
label
$t_j$ on the string which starts from the point $j$ at the bottom. Figure
\ref{braid}
 shows this labelling for the braid
$$\beta=\sigma_1\sigma_2^{-1}\sigma_1\sigma_2^{-1}\sigma_1\sigma_2^{-1}
\sigma_3\in B_4.$$ 

\begin{figure}
\Braid
\caption{The labelled braid $\sigma_1\sigma_2^{-1}\sigma_1\sigma_2^{-1}\sigma_1\sigma_2^{-1}
\sigma_3$.}
\label{braid}
\end{figure}

Write $\overline{C}_i(a)$ for the
$(n-1)\times(n-1)$ matrix  which differs from the unit matrix only in the three
places shown on row $i$, for $1\le i\le n-1$.
\[ \overline{C}_i(a)=\pmatrix{\begin{array}{ccccccc}1&&&&&&\\ &\ddots&&&&&\\
&&1&&&&\\ &&a&-a&1&&\\ &&&&1&&\\ &&&&&\ddots&\\&&&&&&1\\
\end{array} }.\] When $i=1$ or $i=n-1$ the matrix is truncated appropriately to
give two non-zero entries in row $i$. 

Now construct the {\it coloured reduced Burau matrix\/}
$\overline{B}_\beta(t_1,\ldots,t_n)$ of the general braid $$\beta=\prod_{r=1}^l
\sigma_{i_r}^{\varepsilon_r}$$ as a product of matrices $\overline{C}_i(a)$, in
which  $a$ is the label of the current undercrossing string.  This gives
$$\overline{B}_\beta(t_1,\ldots,t_n)=\prod_{r=1}^l
(\overline{C}_{i_r}(a_r))^{\varepsilon_r},$$ where $a_r$ is the label of the
undercrossing string at crossing $r$, counted from the top of the braid.

In the example  shown, where $\beta=\sigma_1\sigma_2^{-1}\sigma_1\sigma_2^{-1}\sigma_1\sigma_2^{-1}
\sigma_3$, the labels $a_1,\ldots,a_7$ are
$t_1,\,t_4,\,t_2,\,t_1,\,t_4,\,t_2,\,  t_4$ respectively and 
$\overline{B}_\beta$ is the
$3\times 3$ matrix product
$$\overline{C}_1(t_1)\overline{C}_2(t_4)^{-1}\overline{C}_1(t_2)
\overline{C}_2(t_1)^{-1}
\overline{C}_1(t_4)\overline{C}_2(t_2)^{-1}\overline{C}_3(t_4).$$

Each braid $\beta$ 
determines a permutation $\pi\in S_n$ by the representation of $B_n$ on $S_n$
in which a string connects position
$j$ at the bottom  to position $\pi(j)$ at the top. In the example
above
$\pi(1)=1,
\pi(2)=2, \pi(3)=4, \pi(4)=3$.

\begin{theorem} The multivariable Alexander polynomial $\Delta_{\hat\beta\cup
A}$, where $A$ is the axis of the closed $n$-braid $\hat\beta$, is given by the
characteristic polynomial
$\det(I-x\overline{B}_\beta(t_1,\ldots,t_n))$ with the identifications
$t_{\pi(j)}=t_j$.
\label{multi}
\end{theorem}

Remarks: (1)\quad Suppose that a link $L$ is presented as the closure of a braid $\beta$ on $n$
strings and that some homomorphism $\varphi:G_L\to H$ is given. Look at the
part of the diagram of $L$ which consists of $\beta$.  The oriented meridians
$x_1,\ldots,x_n$ for the strings at the bottom of $\beta$ determine elements
$\varphi(x_j)\in H$. At each point further up the braid the meridian of a
string which starts at the bottom as string $j$ will  be mapped to the same
element in $H$. Furthermore, when the braid is closed to form $L$ the strings
$j$ and $\pi(j) $ are identified, so that $\varphi(x_j)=\varphi(x_{\pi(j)})$.

When $\pi$ is the product of $k$ disjoint cycles then the link
$L=\hat\beta$ has $k$ components. The variables in the  resulting polynomial
are $x$ and a $k$-element subset of $t_1,\ldots,t_n$ after the identifications
have been made. In the case
$k=1$ the substitution $t_1=\cdots=t_n=t$ in the matrix $\overline{B}_\beta$ gives
the standard reduced Burau matrix $\overline {B}(t)$ for $\beta$.

(2)\quad The discussion of the Alexander polynomials of a link with $k$
components can be done most uniformly in terms of the {\it Alexander
invariant\/} $D_L$ of the link, defined  by \[D_L=\cases {\Delta_L & for
$k>1$,\cr \ds  {\Delta_L(t)\over 1-t}& for $k=1$.}\] The Torres-Fox formula
gives the Alexander invariant of a sublink $L$ of a link $L\cup C$ in terms of
the invariant for $L\cup C$. It says that \[D_L({\bf t}) ={\Delta_{L\cup
C}({\bf t},1) \over 1-\varphi(c)},\] where the meridian of the curve $C$ to be
suppressed is replaced by $1$ and $\varphi(c)$ is the element represented by the
curve $C$ in the complement of $L$, abelianised appropriately. Thus we can
calculate the Alexander invariant of $L=\hat\beta$ from  theorem \ref{multi} by
suppressing the axis and applying the Torres-Fox formula. Put $x=1$ and note
that
$\varphi(A)=t_1t_2\cdots t_n $ in the complement of $L$, to get
\[D_L={\det(I-\overline{B}_\beta(t_1,\ldots,t_n))\over 1-t_1t_2\cdots t_n},\]
making any identifications $t_{\pi(i)}=t_i$ required where strings of $\beta$
belong to the same component of $L$. In the case when $L$ has one component
this gives the well-known formula $\ds
\Delta_L(t)={\det(I-\overline{B}_\beta(t))\over 1-t^n}(1-t)$ quoted earlier.

 \begin{proof}[Proof of Theorem 1] Apply Fox's free differential
calculus  to a presentation of the fundamental group of the closed
braid and axis, very much as in \cite{Birman} or \cite{Morton}. The heart of
the proof lies in relating the fundamental group of the $n$-punctured disk
spanning the axis $A$ which meets $\hat\beta$ at the bottom of $\beta$ to the
corresponding group for the disk at the top of $\beta$. Write $x_1,\ldots,x_n$
for the generators of the group at the bottom represented by meridian loops
around the punctures, and $X_1,\ldots,X_n$ for the meridian loops around the
punctures at the top. It is well known that $X_1,\ldots,X_n$ can be expressed
in terms of
$x_1,\ldots,x_n$ as $X_i=F_\beta(x_i)$, where $F_\beta:F_n\to F_n$ is an
automorphism of the free group $F_n$ determined by the braid $\beta$.
Furthermore the map $\beta\mapsto F_\beta$ from the braid group $B_n$ to
$\hbox{Aut\,} F_n$ is itself a group homomorphism. Then $F_\beta$ is  the
composite of elementary automorphisms corresponding to the elementary braids
making up $\beta$, which are given explicitly by $F_{\sigma_i}(x_i)=x_{i+1}$,
$F_{\sigma_i}(x_{i+1})=x_{i+1}x_i x_{i+1}^{-1}$ and $F_{\sigma_i}(x_j)=x_j,\,
j\ne i,i+1$.

The group of the link $\hat\beta \cup A$ is presented by a generator $x$,
arising from a meridian of $A$, and generators $x_1,\ldots,x_n$ as above, with
$n$ relations $F_\beta(x_i)=x^{-1}x_ix$. The reduced Burau matrix shows up
most naturally by using a different system of generators $g_1,\ldots,g_n$ for
$F_n$, defined recursively by $g_1=x_1,\; g_{i+1}=x_{i+1}g_i$. The element
$g_i$ can be represented by a loop in the disk which encircles the first $i$
punctures as indicated in  figure \ref{disc}.  The automorphism $F_{\sigma_i}$
then satisfies
$F_{\sigma_i}(g_i)=
F_{\sigma_i}(x_ig_{i-1})=x_{i+1}g_{i-1}=g_{i+1}g_i^{-1}g_{i-1}$ and
$F_{\sigma_i}(g_j)=g_j,\; j\ne i$.

\begin{figure}
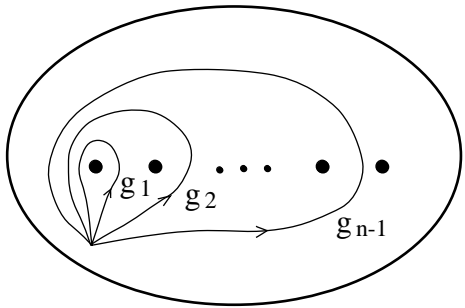

\Disc
\caption{The generators in the punctured disk.}
\label{disc}
\end{figure}

The standard method for finding the Alexander invariant of a link (or knot) $L$
evaluated in
${\bf Z}[H]$ using a homomorphism $\varphi:G_L\to H$ starting from a
presentation of
$G_L$ by
$n+1$ generators $\{g_j\}$ and $n$ relations $\{r_i\}$ is the following.
Calculate the
$n\times(n+1)$ matrix $\ds {\partial r_i\over \partial g_j}$ of free
derivatives, and evaluate the entries in ${\bf Z}[H]$ by applying $\varphi$.
Then delete one column corresponding to a generator
$c$, say, with
$\varphi(c)\ne 1$, and divide the determinant of the remaining matrix by
$1-\varphi(c)$. The result is the evaluation of the Alexander
invariant, that is, the multi-variable Alexander polynomial when $L$ has more
than one component, or 
$\Delta_L(t)/(1-t)$ when $L$ is a knot.
Where a relation is written in the form $r=s$ the entries $\ds
{\partial(r-s)\over\partial g_j}$ work equally well in the matrix.

 The presentation of $G_L$ by generators $g_1,\ldots,g_n$ and $x$ above will
give an immediate calculation of the Alexander invariant, once we know the
entries
$\ds {\partial F_\beta(g_i)\over \partial g_j}$.

The chain rule for free derivatives \cite{Birman} shows that if we write
$\beta=\beta_1\beta_2$ and set $G_i=F_{\beta_2}(g_i)$ then
$F_\beta(g_i)=F_{\beta_1}(G_i)$ and \[ {\partial F_\beta(g_i)\over\partial
g_j}=\sum_{k=1}^n {\partial F_{\beta_1} (G_i)\over\partial G_k}{\partial
G_k\over\partial g_j}.\] We can then compute the matrix  $\ds {\partial
F_\beta(g_i)\over \partial g_j}$ as the product of matrices  $\ds {\partial
F_\gamma(G_k)\over \partial G_j}$, where $\gamma$ runs through the elementary
braids $\sigma_i^{\pm1}$ in $\beta$.

Now when $\gamma=\sigma_i$ we have \[
\begin{array}{rl}F_\gamma(G_i)&=G_{i+1}G_i^{-1}G_{i-1}\\ F_\gamma(G_j)&=G_j,\;
j\ne i.\\
\end{array}
\]

Then \[
\begin{array}{rl}\ds{\partial F_\gamma(G_i)\over \partial G_{i-1}}&=
G_{i+1}G_i^{-1}\\
\ds{\partial F_\gamma(G_i)\over \partial G_i}&=-G_{i+1}G_i^{-1}\\
\ds{\partial F_\gamma(G_i)\over \partial G_{i+1}}&=1\\
\ds{\partial F_\gamma(G_k)\over \partial G_j}&=\delta_{j\,k},
\;\hbox{otherwise}.\\
\end{array}\]

Apply $\varphi$ to the matrix to get
\[\varphi\left({\partial F_\gamma(G_k)\over \partial
G_j}\right)=\pmatrix{\begin{array}{ccccccc}1&&&&&&\\ &\ddots&&&&&\\ &&1&&&&\\
&&a&-a&1&&\\ &&&&1&&\\ &&&&&\ddots&\\&&&&&&1\\
\end{array} },
\] where $a=\varphi(G_{i+1}G_i^{-1})$. This is the value of $\varphi$ on the
meridian at the undercrossing in the current elementary braid
$\gamma=\sigma_i$. Similarly, when $\gamma=\sigma_i^{-1}$ we get the matrix
$$\pmatrix{\begin{array}{ccccccc}1&&&&&&\\ &\ddots&&&&&\\ &&1&&&&\\
&&1&-a^{-1}&a^{-1}&&\\ &&&&1&&\\ &&&&&\ddots&\\&&&&&&1\\
\end{array} }$$ with $a=\varphi(G_iG_{i-1}^{-1})$ which is again the value of
$\varphi$  on the meridian of the undercrossing string.

The generator $g_n$ is clearly unchanged by any automorphism $F_\beta$, so the
last row in each of the matrices of partial derivatives is $(0\cdots 0\;1)$.
The leading $(n-1)\times(n-1)$ submatrix for the partial derivatives of
$F_\beta(g_i)$ is then the product of the leading submatrices from the
constituent elementary braids, and so is the reduced colured Burau matrix of
$\beta$ introduced earlier, with $t_j=\varphi(x_j)$. Thus \[
\varphi\left({\partial F_\beta(g_i)\over\partial
g_j}\right)=\pmatrix{\overline{B}_\beta(t_1,\ldots,t_n)&{\bf v}\cr{\bf
0}&1\cr}=\tilde{B}_\beta(t_1,\ldots,t_n),\;\hbox{say},\] for some column $\bf
v$.  The matrix of free derivatives arising from the presentation of the group
by generators  $g_1,\ldots,g_n$ and $x$ given above is the $n\times(n+1)$ 
matrix whose first $n$ columns, corresponding to the generators
$g_1,\ldots,g_n$, are $\varphi({\partial\over\partial
g_j}(F_\beta(g_i)-x^{-1}g_i x))$. These columns form the matrix
$\tilde{B}_\beta(t_1,\ldots,t_n)-x^{-1}I_n$. They arise from the full matrix of
free derivatives  by deleting the column corresponding to the generator $x$, so
the Alexander invariant for the link is the determinant of this matrix divided
by
$1-x$. Now clearly
$\det(\tilde{B}_\beta(t_1,\ldots,t_n)-x^{-1}I_n)=x^{-n}\det(x\overline{B}_\beta(t_1,\ldots,t_n)-I)(1-x)$.
Then $\det(I-x\overline{B}_\beta(t_1,\ldots,t_n))$, making the identifications
of variables
$t_j=t_{\pi(j)}$ forced by matching the strings at the top and bottom, is the
Alexander polynomial of the closed braid and its axis, up to a power of $x$. 
\end{proof}

The Alexander polynomial of the closed braid itself can then be found from the
Torres-Fox formula as shown above.

\end{document}